\documentclass[12pt]{amsart}
\usepackage{amsmath,amssymb,a4wide}
\usepackage[all]{xy}
\theoremstyle{plain}
\newtheorem{propn}{Proposition}[section]
\newtheorem{thm}[propn]{Theorem}
\newtheorem{lemma}[propn]{Lemma}
\newtheorem{cor}[propn]{Corollary}

\theoremstyle{definition}

\theoremstyle{remark}
\newtheorem{rem}{Remark}
\newtheorem*{ex}{Example}

\def \bnabla{\overline{\nabla}}
\def \V{\mathcal{V}}
\def \H{\mathcal{H}}
\def \bbV{\mathbb{V}}
\numberwithin{equation}{section}

\newcommand{\g}{\mathfrak{g}}
\newcommand{\h}{\mathfrak{h}}

\newcommand{\hs}{\mathfrak{h}_{sym}}
\newcommand{\gs}{\mathfrak{g}_{sym}}
\begin{document}

\title[]{Complex Riemannian foliations of open K\"ahler manifolds}
\date{\today}
\author[T.Murphy]{Thomas Murphy}
\author[P.-A.Nagy]{Paul-Andi Nagy}
\address[T.Murphy]{Department of Mathematics, California State University Fullerton, 800 N. State College Blvd, Fullerton CA 92831, U.S.A. }
\email{tmurphy@fullerton.edu}
\address[P.-A. Nagy]{Department of Mathematics,The University of Murcia, Campus de Espinardo, E-30100 Espinardo,
Murcia, Spain}
\email{paulandi.nagy@um.es}

\begin{abstract}
Classification results for complex Riemannian foliations are obtained.   For open subsets of irreducible Hermitian symmetric spaces of compact type, where one has explicit control over the curvature tensor, 
we completely classify such foliations  by studying the infinitesimal model associated to the canonical connection.  We also establish results for symmetric spaces of non-compact type and  a  general rigidity result for any irreducible K\"ahler manifold.  
\end{abstract}

\maketitle

\section{Introduction}

A central theme in the study of foliations is the classification of Riemannian foliations whose leaves satisfy natural geometric conditions. This question can be traced back to Cartan, who asked if every hypersurface in a round sphere with constant principal curvatures comes from a cohomogeneity one action. Known classifications, such as the work of  Gromoll-Grove \cite{gg} or Kollross \cite{kollross}  have mainly concentrated on two cases: classifying Riemannian foliations of space forms and classifying isometric group actions on Riemannian symmetric spaces under additional hypothesis.

In contrast, our focus is on Riemannian foliations of locally irreducible K\"ahler manifolds  $(M^{2m},g,J), m \geq 2$ where the  leaves are 
complex  submanifolds (a \emph{complex Riemannian foliation}). Such foliations arise in many different contexts, two natural examples arising as 
twistor spaces of 
quaternionic-K\"ahler manifolds and in the classification of closed nearly 
K\"ahler manifolds. Our emphasis is on the local case, where 
existing classification results fail. This paper investigates 
questions which naturally arose from our previous independent 
work \cite{m}, \cite{murphy2}, \cite{Nagy1}, \cite{Nagy3}.

The prototypical example of a twistor fibration comes from the two Hopf fibrations described by the diagram

\vspace{.3in}
$\xymatrix@1{ & & &
\ar[d] S^1 \ar[r] & \mathbb{S}^{4n+3} \ar[r]\ar[rd]& \mathbb{C}P^{2n+1}\ar[d]^{\pi} &\\
& & &  SU(2) \ar[ur] &   &\mathbb{H}P^n } $
\vspace{.3in}

From the inclusion $S^1\rightarrow SU(2)$ the existence of a natural Riemannian submersion $$\pi: \mathbb{C}P^{2n+1} \rightarrow \mathbb{H}P^n$$ with totally 
geodesic leaves $SU(2)/S^1 = \mathbb{C}P^1$ is apparent. Here both 
spaces have their canonical symmetric space metrics. The leaves of a Riemannian submersion are well known to form a Riemannian foliation of the total space. This is the twistor fibration of  $\mathbb{H}P^n$.

A second example of a complex Riemannian foliation discussed in \cite{barbosa}, \cite{besse} (Remark 9.91) is the twistor space of $\mathbb{S}^{2n}$, given by the fibration
$$
H_{n-1} \rightarrow H_n \rightarrow \mathbb{S}^{2n}
$$
where $H_n = SO_{2n+1}/U_n$ and $n>2$. Note that $H_n$ can be identified with   $SO_{2n+2}/U_{n+1}$ equipped with its symmetric space metric, and with this metric the fibration becomes a Riemannian submersion  onto the standard round sphere with complex leaves. This is the  the twistor fibration of $\mathbb{S}^{2n}$. 

Returning to the general case of a complex Riemannian foliation $\mathcal{F}$ on a K\"ahler manifold $M$, denote by $\V$ the leaf tangent distribution and recall that $\mathcal{F}$ is called polar if the distribution 
$\H$ (the orthogonal complement of $\V$ in $TM$) is integrable, and thus totally 
geodesic. Trivial complex foliations are those where every point is a leaf, or where the foliation has just one leaf. Without loss of generality, we will assume that our foliations are non-trivial throughout. 

The main aim of this paper is to describe complex Riemannian foliations of Hermitian locally symmetric spaces. The first step in that direction is the following general 
structure result.

\begin{thm}\label{r2}
For any complex Riemannian foliation $\mathcal{F}$ of a
K\"ahler manifold $M$, let $\mathcal{V}$ denote the distribution tangent to the leaves.
 Then  either $\mathcal{F}$ is totally geodesic, or there is a non-zero integrable
subdistribution  $\mathcal{V}_o\subset \mathcal{V}$, defined on an open set
$\mathcal{O}$ of $M$, whose integral manifolds form a polar complex Riemannian
foliation $\mathcal{F}_0$ of $\mathcal{O}$.
 \end{thm}

This result is optimal: see the Example in Section \ref{red-sect}, where we outline  an algorithm to construct complex, non-polar and non-totally geodesic Riemannian foliations of certain K\"ahler manifolds.
We emphasize that there is no compactness assumption on the ambient manifold, so this parallels existing results for closed K\"ahler manifolds  \cite{Nagy1}. There the proof hinged on a transversal version of the $\overline{\partial}$-lemma for vector bundle valued one-forms. We have 
to develop different arguments to deal with the open case. We also note that our proof does not recover the main result of 
\cite{Nagy1} if $M$ is closed. 

Our main result is the following structure theorem:
\begin{thm}\label{t2}
Let  $M$  be an  open subset of an irreducible Hermitian symmetric space $N$ and $\mathcal{F}$ a complex Riemannian foliation on $M$.  
\begin{itemize}
\item[(i)]  If $N$ has non-negative sectional curvature, then $\mathcal{F}$ is either the twistor fibration of $\mathbb{H}P^n$   restricted to $M\subset \mathbb{C}P^{2n+1}$, or the twistor fibration of $\mathbb{S}^{2n}$ restricted to $M\subset SO_{2n+1}/U_n$.
\item[(ii)] If $N$ has non-positive sectional curvature, then $\mathcal{F}$ is polar.
\end{itemize}
\end{thm}
Theorem \ref{t2} is proved in the following way. Firstly, we use Theorem \ref{r2} 
and some results from almost-K\"ahler geometry to establish that $\mathcal{F}$ must be either polar or totally geodesic. To deal with the latter case we observe that the canonical connection associated with $\mathcal{F}$ must be of Ambrose-Singer 
type. We then analyze the precise relation between the infinitesimal 
model of this connection and the one constructed using the Levi-Civita connection of the ambient space. Once this is established, by direct comparison, the result follows by using a classical result from \cite{on}.

This is a purely local result, whereas all previous classification results assume completeness.  Let us give some applications of the theorem.    Firstly, 
it shows that the canonical twistor
fibration restricted to some open subset $M\subset
\mathbb{C}P^{2n+1}$ provides the only instance when a twistor space of a  positive quaternionic-K\"ahler manifold can be
locally symmetric\footnote{In \cite{chowyang}, it is claimed that the twistor space of any Wolf space is symmetric. Unfortunately this is incorrect:  
the twistor space of $G_2(\mathbb{C}^{m+2})$ is a counterexample. Our result shows it can never happen, even locally. This means the classification of quaternionic-K\"ahler manifolds with non-negative sectional curvature remains an open problem.}.
Combined with recent work in \cite{disc}, Theorem \ref{t2} also completes the classification of complex Riemannian foliations of locally Hermitian symmetric spaces with non-negative sectional curvature. Theorem \ref{t2} strengthens  results
in \cite{m}, where this result was established using different
techniques for complex Riemannian foliations of open subsets of
$\mathbb{C}P^n$ with maximal dimension.

\begin{rem} It is natural in light of Theorem 1.1(ii) to ask whether there are any explicit examples of complex polar Riemannian foliations of a Hermitian symmetric space on non-compact type. Interestingly there are many known examples of polar homogeneous foliations of such manifolds (see \cite{bd}), but all these examples are known to not be complex (see Proposition 4.3 of \cite{bd}). We hope to return to this question in a future work.
\end{rem}

\bigskip

\noindent\textbf{Acknowledgements}  T.M. was supported by an A.R.C. grant 
whilst at the Universite Libre de Bruxelles, a Britton Fellowship 
whilst at McMaster University, and by a startup research grant 
from California State University Fullerton. The research 
of P.-A.N. has been partially supported by the 
grant MTM2012-34037(MICINN). We also wish to thank R. Bryant,  A. Swann, 
McKenzie Y.K. Wang and D. Yang for helpful conversations. 
We thank the referee for many helpful comments on 
our paper.  We especially thank S. Salamon for 
pointing out we missed the twistor space of $\mathbb{S}^{2n}$ in the first draft of this work.

%Although this is not strictly speaking necessary for the arguments in this work we %have included an Appendix to this paper. It containes a general, very short %criterium for an infinitesimal model to be regular, of a type we could not find in %the litterature.
\section{Preliminaries} \label{prel}
Let $(M^{2m},g,J), m \geq 2$ be a connected K\"ahler manifold. Throughout this paper the metric $g$ will be assumed to be locally irreducible. That is the connected 
component $Hol^0(g)$ of its holonomy group acts irreducibly on some, and therefore all, tangent spaces. In particular, this happens in case $M$ is an open subset of a simply connected, irreducible, Hermitian symmetric space.
%For, any $Hol^0(g)$-invariant splitting of $TM$ at some point extends, via %parallel 
%transport, to a 
%Since this is a local concept, any open subset of an irreducible Hermitian %symmetric space will have this property. \footnote{more detail?}

Now assume that $M$ is equipped with a foliation $\mathcal{F}$, with tangent leaf distribution $\V$. In what follows $\mathcal{F}$ will be assumed to be {\it{complex}}, that is $J\V=\V$. Let $\H$ denote the orthogonal complement with respect to $g$ of $\V$ in $TM$; it clearly satisfies $J\H=\H$. The foliation $\mathcal{F}$ is {\textit{Riemannian}} if it respects the metric: that is the distance between leaves is constant. In terms of the Lie derivative of the metric $g$ this can be expressed as 
\begin{equation*}
(L_{V}g)(X,Y) = 0
\end{equation*} whenever $V$ respectively $X,Y$ are smooth sections of $\V$, respectively $\H$. In any K\"ahler manifold there are two trivial examples of complex Riemannian foliations; the foliation with only one leaf, and the foliation
where each leaf is a point. Equivalently $\V=TM$ or $\V=\{0\}$. We will 
assume our foliations are non-trivial throughout. 
 
Assuming henceforth that $(M^{2m},g,J)$ is equipped with a complex Riemannian foliation 
$\mathcal{F}$, we record a few basic facts we will need in what follows. 
Throughout this paper we will denote by $V,W,$ sections of $\mathcal{V}$   
and by $X,Y,Z$ etc. sections of $\mathcal{H}$. Generic sections of $TM$ will be denoted by $U_1,U_2$. Let $\nabla$ be the Levi-Civita connection
of $g$. Then 
$$ \overline{\nabla}=2\nabla+\pi^{\V} \circ \nabla \pi^{\V}+\pi^{\H} \circ \nabla \pi^{\H}
%(\nabla_EF_{\mathcal{V}})_{\mathcal{V}}+(\nabla_EF_{\mathcal{H}})_{\mathcal{H}}
$$
defines a metric connection, that is $\bnabla g=0$ on $M$.  Here $\pi^{\V}:TM \to \V$ (resp., $\pi^{\H}:TM \to \H$) are the orthogonal projection maps onto $\V$ 
(resp., $\H$). Because $\nabla J=0$ and since $\V$ and $\H$ are both $J$-invariant it follows that $\bnabla$ is also Hermitian, that is $\bnabla J=0$. 

The main feature of the connection $\bnabla$ is to preserve the distributions $\mathcal{V}$ and $\mathcal{H}$. Note that $\bnabla$ has torsion, from which 
one can interpret the assumptions made about the foliation in the following manner. If $T$ and $A$ are the O'Neill tensors of the foliation then $\nabla$ and $\overline{\nabla}$ are related by
\begin{equation} \label{comp-co}
 \begin{array}{lr}
\nabla_XY=\overline{\nabla}_XY+A_XY, \ \nabla_XV=\overline{\nabla}_XV+A_XV \\
\nabla_VX=\overline{\nabla}_VX+T_VX, \ \nabla_VW=\overline{\nabla}_VW+T_VW.
\end{array}
\end{equation}
Recall that $A$ is skew-symmetric
on $\mathcal{H}$ since the foliation is Riemannian whereas $T$ is symmetric on $\mathcal{V}$, since the latter is integrable. Assuming that $\bnabla$ is Hermitian yields
information about the complex type of the tensors $A$ and $T$, namely
\begin{equation} \label{21}
\begin{split}
A_X(JY)=J(A_XY), & \ A_{JX}V=-J(A_XV)=-A_X(JV) \\
T_{JV}W=J(T_VW), & \ T_{JV}X=-J(T_VX)=-T_V(JX).
\end{split}
\end{equation}
We also have $A_{JX}JY=-A_XY$ and $T_{JV}JW=-T_VW$.  \par

We will use now the K\"ahler structure on $M$, together with suitable
curvature identities to collect some geometric information about the tensors $A$ and $T$. Of central importance is the following result.
\begin{propn} \label{l1} %been upgraded to Proposition
Let $X,Y,Z$ be in $\mathcal{H}$ and $V,W$ in $\mathcal{V}$. Then
\begin{itemize}
\item[(i)] $(\overline{\nabla}_XA)(Y,Z)=0$
\item[(ii)] $\langle A_XY, T_VZ \rangle=0$
\item[(iii)] $\langle (\overline{\nabla}_VA)(X,Y),W \rangle=\langle (\overline{\nabla}_WA)(X,Y),V \rangle$
\item[(iv)] $ 2\langle (\overline{\nabla}_VA)(X,Y), W \rangle=\langle(\overline{\nabla}_YT)(V,W),X \rangle-\langle(\overline{\nabla}_XT)(V,W), Y \rangle$
\item[(v)] $\overline{\nabla}_{A_XY}A=0$.
\end{itemize}
\end{propn}
\proof
Parts (i) to (iv) have been proved in \cite{Nagy3}. To prove (v) we differentiate (ii) with respect to
$\overline{\nabla}$, in direction of $Z_1$ in $\mathcal{H}$. Then, by taking (i) into account
$$\langle A_XY, (\overline{\nabla}_{Z_1}T)_VZ\rangle =0 .$$ 
Now applying (iv), skew symmetrisation in $Z_1$ and $Z$ yields $$\langle A_XY, (\overline{\nabla}_VA)_{Z}Z_1\rangle=0.$$ The claim follows now from the
symmetry property for $\overline{\nabla} A$ in (iii).
\endproof
Let $R$ be the curvature tensor of the connection $\nabla$ defined according to  the convention $R(U_1,U_2)=\nabla^2_{U_2,U_1}-\nabla^2_{U_1,U_2}$.
As a consequence of Proposition \ref{l1}, some of the components of the Riemann curvature tensor become algebraic in $A$ and $T$, as showed below.  \par
\begin{cor} \label{curv-s}We have the following. 
\begin{align*}
R(X,Y,Z,V)= &0, \\
R(V,W,X,Y) = & \langle A_XV, A_YW \rangle-\langle A_XW, A_YV \rangle- \langle T_VX,T_WY \rangle+\langle T_WX,T_VY \rangle,\\
R(V,X,W,Y) = & \langle (\overline{\nabla}_XT)(V,W),Y\rangle+\langle (\overline{\nabla}_VA)(X,Y),W \rangle  -\langle T_VX,T_WY\rangle\\
& +\langle A_XV, A_YW \rangle.
\end{align*}
\end{cor}
\proof These follow from Proposition \ref{l1} and  the O'Neill formulae in \cite{besse} applied
to curvature terms of the type listed above.
\endproof
To finish this section we set a few notational conventions. Whenever 
$W$ is a vector space equipped with an inner product $g$ we denote by 
$\mathfrak{so}(W)$ the space of endomorphisms of $W$, skew-symmetric with respect to 
$g$. We will systematically identify $\mathfrak{so}(W)$ and $\Lambda^2W$ via 
$F \mapsto g(F \cdot, \cdot)$. If moreover $J$ is a linear complex structure on $W$, compatible with $g$, we indicate with $\mathfrak{u}(W)$ the sub-space of 
Hermitian endomorphisms in $\mathfrak{so}(W)$. No reference to $g$ nor $J$ will be made, unless otherwise specified.

\section{Reduction results} \label{red-sect}
In this section $(M^{2m},g,J),m \geq 2$ denotes a connected K\"ahler manifold whose tangent bundle splits as $TM=\V \oplus \H$ as considered in the previous section. We will first prove Theorem \ref{r2}, based on a geometric interpretation of Proposition \ref{l1}. Let the symmetric tensor $r:\V \to \V$ be defined by 
$r(V):=-\sum \limits_{i}A_{e_i}A_{e_i}V$ whenever $\{e_i\}$ is a local orthonormal 
frame in $\H$. Note that $r$ is positive and $J$-invariant. 

At each point in $M$ define the $J$-invariant subspaces 
%Consider the locally \footnote{those are global!!} defined distributions
$$ \mathcal{V}_1:=A_{\mathcal{H}}\mathcal{H} \text{, and } \mathcal{V}_0:=\mathcal{V}_1^{\perp}\subset \mathcal{V}.
$$
It is easy to see that $\ker(r)=\V_0$, hence $\V_1=Im(r)$. 
Note that the rank of $r$ may not be constant over $M$. However, one can choose an open dense set $\mathcal{O}$ in $M$ such that $r$ has constant rank over each connected component of $\mathcal{O}$. Over each such open set, 
$\V_0$ and $\V_1$ are thus genuine complex distributions.

\begin{propn} \label{polar}
Let $(M^{2m},g,J),m \geq 2$ be a K\"ahler manifold equipped with a
complex Riemannian foliation $\mathcal{F}$. Then either $\mathcal{F}$ is totally geodesic, or there is an open dense 
subset $\mathcal{O} \subset M$ such that ${\mathcal{V}}_0$ induces a polar Riemannian complex foliation over each connected component of $\mathcal{O}$.
\end{propn}
\proof
Working locally as outlined above, we prove first integrability for the distribution $\mathcal{V}_0$. Choose sections $V_0, W_0$ in $\mathcal{V}_0$. By definition $A_\mathcal{H}\mathcal{H}$ and $\mathcal{V}_0$ are orthogonal, and so
$$\langle (\overline{\nabla}_{V_0}A)_XY,W_0\rangle=\langle \overline{\nabla}_{V_0}(A_XY),W_0 \rangle=-
\langle A_XY, \overline{\nabla}_{V_0}W_0 \rangle.$$
Since the left hand side above is symmetric in $V_0$ and $W_0$ by Proposition \ref{l1} (iii) the integrability of $\mathcal{V}_0$ follows. The next step is to prove that $\mathcal{V}_1\oplus \mathcal{H}$ is integrable and totally geodesic.  For $X,Y$ in $\mathcal{H}$ we have $\nabla_XY=\overline{\nabla}_XY+A_XY$ in $\mathcal{H} \oplus  \mathcal{V}_1$. 
Proposition \ref{l1} (ii) shows that
$$T(\mathcal{V}_1,\mathcal{V})=0.$$
 At the same time, part (v) in the same Proposition guarantees that $\overline{\nabla}_{V_1}A=0$ for $V_1$ in $\mathcal{V}_1$. In particular $\mathcal{V}_1$ is totally geodesic with respect to $\overline{\nabla}$. Therefore $$\nabla_{V_1}W_1=\overline{\nabla}_{V_1}W_1 \in
\mathcal{V}_1.$$

If $(X,V_1)$ is in $\mathcal{H} \times \mathcal{V}_1$ we have $\nabla_{V_1}X=\overline{\nabla}_{V_1}X$ in $\mathcal{H}$ again by using $T(\mathcal{V}_1, \cdot)=0$. Finally, note that $\mathcal{V}_1$ is parallel in direction of $\mathcal{H}$, with respect to $
\overline{\nabla}$. This is a consequence of Proposition \ref{l1} (i). Thus for a pair $(V_1,X)$ as above $\nabla_XV_1=\overline{\nabla}_XV_1+A_XV_1$ belongs
to $\mathcal{V}_1 \oplus \mathcal{H}$.The claim is now fully proved.

\endproof

This establishes Theorem \ref{r2}, because if $\mathcal{V}_0$ vanishes then the foliation is totally geodesic. The result parallels  the situation for closed K\"ahler manifolds, where $\mathcal{F}$ must be either totally geodesic or polar. 
The natural question now is whether the above result could be improved to show the second author's result also holds in the non-compact case. Certainly the techniques in \cite{Nagy1} will not work here. In fact the above theorem is optimal, as we will 
demonstrate by showing how to construct a wide family of examples  $(M,\mathcal{F})$ with $\mathcal{F}$ a complex Riemannian foliation that is neither polar nor totally geodesic. To this end, we follow a construction in \cite{nagychiossi}. 

\begin{ex}Let $N$ be an arbitrary K\"ahler manifold and $\mathbb{U}$ the unit ball in $\mathbb{C}$. Take an arbitrary holomorphic function $f: N\rightarrow \mathbb{U}$. Set
$$
\Phi = \bigg( \begin{array}{cc}  Re(f) \  Im(f) \\ 
Im(f) \ - Re(f) \end{array}\bigg).
$$
Then $\Phi\in S^2(T\mathbb{U})$ with respect to the standard basis. Now consider the product bundle  with projection map $\pi: M=
\mathbb{U}\times N \rightarrow N$. Fix a  K\"ahler metric $g_N$ on $N$, and choose the standard complex structure on $T\mathbb{U}$, called $J_0$. Take the metric
$$
g = g_0\bigg((1+\Phi)^{-1}(1-\Phi)\cdot, \cdot \bigg) + g_N
$$
together with the compatible complex structure
$$
J = (1-\Phi)^{-1}J_0(1+\Phi)  + J_N,
$$
Then $(M,g,J)$ is K\"ahler, and   $$TM= \mathbb{C}\oplus TN = \mathcal{V}_0 \oplus
\mathcal{H}_0$$ is a polar foliation. The key point is that there is no
restriction on $N$(apart from it being non-compact), so we can choose it to be a K\"ahler manifold
equipped with a totally geodesic Riemannian foliation, inducing a
splitting $\mathcal{H}_0 = TN = \mathcal{V}_1 \oplus \mathcal{H}$. Taking now the splitting
$\mathcal{V}\oplus \mathcal{H}$ of $T(\mathbb{U}\times N)$,  with $\mathcal{V}$ equal to $\mathcal{V}_0 \oplus
\mathcal{V}_1$,  it is easy to check that $\mathcal{V}$ is a complex Riemannian
foliation which is not polar unless $\mathcal{V}_1$ splits off. Of course one also sees that $\mathcal{V}_0$
induces a polar foliation as the above theorem states. Concrete examples are given by taking $N$ to be an open part of the twistor space of a quaternionic K\"ahler manifold with positive scalar curvature. 
\end{ex}
 Our next result reduces the study  of complex Riemannian foliations to  a special case for open subsets of Hermitian symmetric spaces of compact type.

\begin{propn} \label{prop1}
Let  $M$  be an open set of an irreducible Hermitian symmetric space and $\mathcal{F}$ a complex Riemannian foliation on $M$. If $M$ has non-negative sectional curvature, then $\mathcal{F}$ is totally geodesic.
\end{propn}
\proof
Consider the splitting $$TM=\mathcal{V}_0 \oplus (\mathcal{V}_1 \oplus \mathcal{H})$$ over $M$. Because both factors are integrable distributions and $(g,J)$ is K\"ahler, it is a well known fact (see e.g. \cite{apo2}) that  
%The integrability ofboth factors implies that  
$(g, \tilde{J})$ is almost-K\"ahler, where
$$ \tilde{J}\vert_{\mathcal{V}_0}=-J, \ \tilde{J}\vert_{\mathcal{V}_1\oplus \mathcal{H}}=J.
$$
That is, $\tilde{J}$ is $g$-orthogonal and the K\"ahler form $\tilde{\omega} = g(\tilde{J} \cdot, \cdot)$ is closed. 
%This follows from the formula
%$$
%d\omega (E_1, E_2, E_3) = E_1\cdot \omega(E_2, E_3) - \omega([E_1, E_2], E_3) - %\omega(E_1, [ E_2, E_3]).
%$$
%Via a  case-by-case analysis of the various possibilities (taking each $E_i$ to %either lie in $\V_0$ or $\mathcal{V}_1 \oplus \mathcal{H}$), one can exploit the %fact $d\omega = 0$ to prove that  $\tilde{\omega}$ is closed. For example, if all %$E_i$ lie in $\V_0$ then
%$$
%d\tilde{\omega}(E_1, E_2,E_3) = - d\omega(E_1, E_2, E_3) = 0.
%$$ 
%One uses the integrability of $\V_0$ to deduce that $J$ changes by a minus sign in %all three components. 
Since  $(M,g)$ is an open part of a symmetric space of compact type its isotropic sectional curvatures must be non-negative.
It follows from  \cite{apo1}, Proposition 1, (iii) that $\nabla \tilde{J}=0$. Equivalently the distributions $\mathcal{V}_0$ and $\mathcal{V}_1 \oplus \mathcal{H}$ are parallel with respect to the Levi-Civita connection of $g$. As the latter is locally irreducible 
either $\mathcal{V}_0=0$, corresponding to the totally geodesic case or $\mathcal{V}_1 \oplus \mathcal{H}=0$, which happens when the foliation is trivial.
\endproof
To finish this section, we will prove the last part in Theorem \ref{t2},relying on 
the following elementary observation. 
\begin{lemma} \label{non-pos}
Let $(M^{2m},g,J)$ be a K\"ahler manifold with non-positive sectional curvature. 
Then any totally geodesic complex Riemannian foliation comes from a local product.
\end{lemma}
\begin{proof}
Consider the splitting $TM=\V \oplus \H$. Because $T=0$ and $\bnabla A=0$ 
we get $R(V,X,V,X)=\Vert A_XV \Vert^2$ by Corollary \ref{curv-s}. It follows that 
$A=0$ as claimed.

\end{proof}
$\\$
{\bf{Proof of Theorem \ref{t2}, (ii)}}\\
We have the splitting
$$
TM = \mathcal{V}_0 \oplus (\mathcal{V}_1 \oplus \mathcal{H}).
$$
Since $M$ has non-positive sectional curvatures and the distribution $\mathcal{V}_1 \oplus \mathcal{H}$ is totally geodesic, the sectional curvatures of any of its integral manifolds are non-positive with
respect to the induced metric. Each integral manifold carries a totally geodesic Riemmanian foliation corresponding to the leaves of $\mathcal{V}_1$ since $T(\mathcal{V}_1, \V) = 0$. The corresponding O'Neill tensor is the restriction 
of $A$ to $\V_1 \oplus \H$. It follows by Lemma \ref{non-pos} that $A$(and thus 
$\V_1=A_{\H}{\H}$) vanish on any integral submanifold of $\V_1 \oplus \H$, thus everywhere in $M$. Since $\V_1=0$ means that the foliation is polar the claim is proved.
\endproof

\begin{rem}
Irreducible compact Hermitian symmetric spaces do not admit non-trivial polar foliations  \cite{lytchak}. 
An alternative argument for the case of $\mathbb{C}P^n$ can be found in  \cite{m}. 
\end{rem}
\section{Proof of Theorem \ref{t2}(i)}

Consider a locally irreducible and connected K\"ahler manifold $(M^{2m},g,J), m\geq 2$, equipped with a totally geodesic, complex, Riemannian foliation with leaf tangent space $\V$. Moreover we assume that $g$ is a locally symmetric metric, that is $\nabla R=0$, which moreover has non-negative sectional curvature.  Since $T=0$ we have $\bnabla A=0$ by Proposition \ref{l1}, (iv).

For an element  $V\in \mathcal{V}$ define $\gamma_V:\mathcal{H} \to \mathcal{H}$ by
$$
\gamma_VX=A_XV.
$$
Then $\gamma_V$ is skew-symmetric and moreover 
$$
\gamma_{JV}= -\gamma_VJ=J\gamma_V.
$$
\begin{rem} \label{dim}
If $\H$ has real rank $2$, the space $\{F \in \mathfrak{so}(\H):FJ+JF=0\}$ 
vanishes, hence so does $\gamma$. This is prevented by the local irreducibility  
assumption on $g$ thus $\dim_{\mathbb{R}}\H \geq 4$. In particular we must have 
$m \geq 3$.
\end{rem} 

\begin{lemma}\label{gammairr} At any given point of $M$ we have that 
\begin{itemize}
\item [(i)] $\gamma : \V \to \Lambda^2\H$ is injective,
\item[(ii)] $\mathcal{V} = A_{\mathcal{H}}\mathcal{H}$. 
\end{itemize}
\end{lemma}
\proof
To show (i), consider the subspace $\mathcal{D} = \lbrace V\in \V : \gamma_V = 0 \rbrace$. Since $\gamma$ is $\overline{\nabla}$ parallel, this is invariant under the holonomy group $Hol(\bnabla)$ of $\bnabla$. It therefore extends, by parallel transport, to a $\bnabla$-parallel distribution, still denoted by $\mathcal{D}$ in $TM$. From \eqref{comp-co}, the intrinsic torsion of $\bnabla$ vanishes on $\mathcal{D}$, which is henceforth $\nabla$-parallel. As $g$ is locally irreducible it follows that $\mathcal{D}=0$. \\
For (ii), since the orthogonal complement of $A_{\mathcal{H}}\mathcal{H}$ coincides with 
$\ker(\gamma)$, the claim follows from (i).
\endproof
\begin{rem} \label{nk}
Define an almost-Hermitian structure on $M$ by 
\begin{equation*}
\begin{split}
&g^{nk}\bigg|_{\mathcal{V}} = \frac{1}{2}g \bigg|_{\mathcal{V}}, \ \ \  g^{nk}\bigg| _{\mathcal{H}} = g \bigg|_{\mathcal{H}}, \ \ \ \  g^{nk}(\mathcal{H},\mathcal{V}) = 0 \\
&J^{nk}=-J_{\vert \V}+J_{\vert \H}.  
\end{split}
\end{equation*}
As it is well known (see \cite{nagyagag})the pair $(g^{nk},J^{nk})$ is nearly-K\"ahler. It can be used to give another interpretation for the connection 
$\bnabla$. If $\nabla^{nk}$ denotes the Levi-Civita connection of $g^{nk}$ 
an elementary computation shows that   
\begin{equation*}
\overline{\nabla} = \nabla^{nk} + \frac{1}{2} \left(\nabla^{nk} J^{nk} \right) J^{nk}.
\end{equation*}
The right hand side above is precisely the canonical Hermitian connection of the nearly-K\"ahler structure. This  observation will be used as a guideline on how to adapt to our situation (see  Propositions \ref{irrep-fn} and \ref{irrep-V} in the next section) some purely representation-theoretical aspects of the classification of nearly-K\"ahler metrics which have been developed in \cite{handbook}.
\end{rem}
\subsection{Initial Computations} \label{basics}
Let us now calculate first-order information on the structure of the curvature tensor of $g$, using the fact that $\nabla R=0$. In what follows $[\cdot, \cdot]$ denotes the commutator on operators. We let $\overline{R}$ be the curvature tensor of the connection $\bnabla$, defined according to $\overline{R}(U_1,U_2)=
\bnabla^2_{U_2,U_1}-\bnabla^2_{U_1,U_2}$ whenever $U_1,U_2 \in TM$. Below 
$R^{\H}(X,Y)$ (resp., $[A_X,A_Y]^{\H}$) denote the component of $R(X,Y)$ 
(resp., $[A_X,A_Y]$) on $\Lambda^2\H$.
\begin{lemma} \label{curv2} The following equations hold:
\begin{itemize}
\item[(i)] $\gamma_{R(V_1,V_2)V_3}=[[\gamma_{V_1},\gamma_{V_2}], \gamma_{V_3}],$
\item[(ii)]  $\gamma_V R^{\H}(X_1,X_2)=[A_{X_1},A_{X_2}]^{\H} \gamma_V-\gamma_{[A_{X_1},A_{X_2}]V},$
\item[(iii)] $\bnabla R=0,$
\item[(iv)] $\bnabla \ \overline{R}=0$.
\end{itemize}
\end{lemma}
\proof
Part (i) is a consequence of the differential Bianchi identity and holds independently of the
assumption $\nabla R=0$. See  \cite{Nagy1} for details. \\
Part (ii) is a consequence of the fact our metric is locally symmetric:
$$(\nabla_{X_4}R)(X_1,X_2,X_3,V)=0.$$ 
Expanding and using the fact that $R(X_1,X_2,X_3,V)=0$ yields
\begin{align*}
0 = & R(A_{X_4}X_1,X_2,X_3,V)+R(X_1,A_{X_4}X_2,X_3,V) \\
&+R(X_1,X_2,A_{X_4}X_3,V) +R(X_1,X_2,X_3,A_{X_4}V).
\end{align*}
Expressing the first three terms above by means of Corollary \ref{curv-s} a short algebraic computation  establishes the second claim. \\
For (iii), differentiate with respect to $\bnabla$ in (ii). Since $\bnabla A=0$ it follows by Lemma \ref{gammairr} that $(\bnabla R)(X,Y)=0$. At the same time, applying the same argument in (i) shows that the component of $\bnabla R$ in $TM \otimes \Lambda^2\V \otimes \Lambda^2\V$ vanishes. The components of $R$ on $(\V \wedge \H) \otimes \Lambda^2M$ and $(\Lambda^2 \V \otimes \Lambda^2\H)
\oplus (\Lambda^2 \H \otimes \Lambda^2\V)$ are explicitly determined in terms of $A$ only (see Corollary \ref{curv-s}). Therefore they must be $\bnabla$-parallel and the claim is proved.\\
(iv) a straightforward computation only based on $\bnabla A=0$ leads to 
\begin{equation} \label{curv11}
\overline{R}(X,Y)=R(X,Y)+[A_X,A_Y]
\end{equation} 
for all $X,Y \in \H$. Differentiating and using (iii) shows that $\bnabla 
\ \overline{R}(X,Y)=0$. Because $\overline{R}$ is symmetric in pairs $\overline{R}(V,X)=0$. Finally $\overline{R}(V,W)=R(V,W)$ as $\V$ is totally geodesic so the claim follows from (iii).
\endproof
In particular $\bnabla$ is an Ambrose-Singer connection since it preserves  
its own torsion and curvature tensors. As an immediate consequence 
$$\h_{nk}=span \{ \overline{R}(U_1,U_2): U_1, U_2 \in \bbV \} $$
is a Lie subalgebra of $\mathfrak{u}(\bbV)$. Here $\bbV:=T_pM$, where $p$ is a fixed point in $M$, is equipped with its standard metric and complex structure. In fact, $\h_{nk} $ coincides with the Lie algebra of $Hol^0(\bnabla)$, the connected component of the holonomy group of $\bnabla$.

In the rest of this section we will 
describe the main properties of the representation $(\h_{nk},\mathbb{V})$. Clearly $\h_{nk}$ preserves $\V$ and $\H$ . Moreover, 
\begin{propn} \label{irrep-fn}
The representation $(\h_{nk},\H)$ is irreducible.
\end{propn}
\begin{proof}
Assume $\mathcal{H} = \mathcal{H}_1\oplus\mathcal{H}_2$ is an orthogonal, $J$-invariant and $\h_{nk}$-invariant splitting. Then 
$ \overline{R}(X,Y,X_1, X_2)=0$
for all $X,Y \in\mathcal{H}$ and  $X_i\in \mathcal{H}_i$. In particular
\begin{align*}
0 =& \overline{R}(X_1,X_2,X_1, X_2)
= R(X_1,X_2,X_1,X_2) + \langle [A_{X_1}, A_{X_2}]X_1, X_2\rangle\\
=& sec_g(X_1,X_2) + \|A_{X_1}X_2\|^2
\end{align*}
by using the comparaison formula \eqref{curv11}. As $g$ has non-negative sectional curvatures, it follows that $A_{\mathcal{H}_1}\mathcal{H}_2 = 0$. From now on the proof proceeds in the spirit of \cite{handbook}, Section 5.2, so the details will only be sketched. Indeed, the algebraic Bianchi identity for $\overline{R}$(obtained for instance by taking the cyclic sum in \eqref{curv11}) leads to 
$\overline{R}(X_1,Y_1,X_2,Y_2)=2\langle A_{X_2}A_{X_1}Y_1, Y_2\rangle $
for all $X_1,Y_1 \in \H_1$ and $X_2,Y_2 \in \H_2$. Since the left hand side is 
$J$-invariant in $(X_2,Y_2)$ it follows easily that the right hand side vanishes, that is the subspaces $\V_1:=A_{\H_1}\H_1$ and $\V_2=A_{\H_2}\H_2$ are orthogonal. 
Moreover $\V=A_{\H}\H=\V_1 \oplus \V_2$ by Lemma \ref{gammairr}. Now observe that 
$\gamma_{\V_1}\H_2$ is orthogonal to $\H_1$ (since $A_{\H_1}\H_2=0$) and also to 
$\H_2$ (since $A_{\H_2}\H_2=\V_2 \perp \V_1$). Therefore $\gamma_{\V_1}\H_2=0$ and 
similarly $\gamma_{\V_2}\H_1=0$. Consider the $\h_{nk}$-invariant splitting 
%The $\h_{nk}$-invariant splitting 
\begin{equation} \label{split-int}
\bbV=(\V_1 \oplus \H_1) \oplus (\V_2 \oplus \H_2).
\end{equation} 
It extends, by parallel transport, to a $\bnabla$-parallel splitting of $TM$ over any simply connected open subset 
$\mathcal{O}$ in $M$ with $p \in \mathcal{O} $, due to 
$Hol(\bnabla_{\vert \mathcal{O}})=Hol^0(\bnabla_{\vert \mathcal{O}})$. At the same time, the algebraic facts proved above amount to saying that \eqref{split-int} is preserved by the intrinsic torsion tensor 
of $\bnabla$, given by $\eta:=\bnabla-\nabla$; that is 
$\eta_U(\V_i \oplus \H_i) \subseteq \V_i \oplus \H_i$ for $i=1,2$ whenever 
$U \in \bbV$. Thus our splitting of $TM$ is also $\nabla$-parallel, which contradicts the local irreducibility of $g$.
\end{proof}
In the same vein we have:
\begin{propn} \label{irrep-V}
The representation $(\h_{nk},\V)$ is irreducible.
\end{propn}
\begin{proof}
Assume that $\V=\V_1 \oplus \V_2$ where $\V_k, k=1,2$ are $\h_{nk}$-invariant, $J$-invariant and orthogonal with respect to $g$. Then 
$\overline{R}(X,Y,V_1,V_2)=0$ for all $X,Y \in \H$ and $V_i \in \V_i$ for $i=1,2$. 
However from \eqref{curv11} and Corollary \ref{curv-s} we have 
$\overline{R}(X,Y,V_1,V_2)=2 \langle [\gamma_{V_1}, \gamma_{V_2}]X,Y\rangle$. Thus 
$[\gamma_{V_1}, \gamma_{V_2}]=0$ and further $\gamma_{V_1}\gamma_{V_2}=0$ by using the $J$-invariance properties of the tensor $\gamma$. Equivalently, the spaces 
$\gamma_{\V_1}\H$ and $\gamma_{\V_2}\H$ are orthogonal within $\H$. Since they 
are also $\h_{nk}$-invariant, the irreducibily of $(\h_{nk},\H)$ ensures that 
one of them must vanish. Then by Lemma \ref{gammairr} we have either $\V_1=0$ or $\V_2=0$ as claimed.
\end{proof}

\subsection{Comparing the infinitesimal models} \label{sect 5}
Because $\bnabla $ is an Ambrose-Singer connection, Nomizu's construction 
allows one to build its associated infinitesimal model. We refer to \cite{tricerri} for the general theory. This a purely algebraic construction, so in what follows we adopt the following notation and conventions. The vector space $\bbV$ will be always equipped with its canonical metric $g$ and linear complex structure $J$. Elements in $\bbV$ 
will be denoted by upper-case letters, $U_1,U_2$ etc., while for elements in $\V$ respectively $\H$ we use lower case letters, for example $v$ respectively $x$. Elements in $\mathfrak{u}(\bbV)$ will generally be denoted by $F,G$ etc.

The infinitesimal model for the Ambrose-Singer connection $\bnabla$ is the Lie algebra 
$$\mathfrak{g}_{nk} = \mathfrak{h}_{nk} \oplus \mathbb{V}
$$
with Lie bracket
\begin{align*}
[U_1,U_2]_{nk} & = \tau(U_1,U_2) + \overline{R}(U_1,U_2)\\
[F,U]_{nk} &= FU,\\
[F,G]_{nk} &= FG - GF,
\end{align*}

\noindent where $U_1,U_2,U\in \mathbb{V}$ and $F,G\in \mathfrak{h}_{nk}$. Here $\tau$ is the torsion tensor of $\bnabla$. For convenience, we recall here that 
$$\tau(x,y)=-2A_xy, \ \tau(v,x)=A_xv.
$$ 

Now consider the infinitesimal model for the locally symmetric space $(M,g)$. Letting $\mathfrak{h}_{sym}:=span \{R(U_1,U_2): U_1,U_2 \in \mathbb{V}\}$ we have 
$\gs=\hs \oplus \mathbb{V}$ with Lie bracket given by 
\begin{align*}
 [U_1,U_2]_{sym} &= R(U_1,U_2), \\
 [F,U]_{sym}&= FU,\\
 [F,G]_{sym}&=FG-GF.
 \end{align*}

The strategy for the proof of Theorem \ref{t2}, part(i), involves thoroughly comparing the 
infinitesimal models $(\g_{nk},\h_{nk})$ and $(\gs,\hs)$. 
\begin{rem} \label{old-line} It is possible to see that the infinitesimal model $(\g_{nk},\bbV)$ is {\textit{regular}}
(see \cite{tricerri} for the definition); then one can conclude that $(M,g^{nk},J^{nk})$ 
arises locally as an open part of a closed nearly-K\"ahler manifold $(G^{nk}\slash H^{nk},\tilde{g}^{nk},\tilde{J}^{nk})$. Performing the canonical variation in direction of $\V$, reversing the construction in Remark \ref{nk}, one obtains a K\"ahler structure $(h,I)$ 
on $G^{nk}\slash H^{nk}$ which is moreover locally symmetric and therefore of type $G^{sym}
/H^{sym}$. Compact, simply connected homogeneous nearly-K\"ahler 
manifolds have been explicitly classified in \cite{CGD}. One could then compare their list to that of Hermitian symmetric spaces of compact type. To do so,
one needs a priori information on how the pairs $(G^{nk},H^{nk})$ and $(G^{sym},H^{sym})$ relate to each other; in particular that $G^{nk}$ is a {\textit{strict}} subgroup of $G^{sym}$ in the sense that $G^{nk} \neq G^{sym}$. After describing the relation between the pairs $(\g_{nk},\h_{nk})$ and $(\g_{sym},\h_{sym})$ in detail below, we have realized a much shorter argument is available, based on a classical result of Onishchik \cite{on}.
\end{rem}
We begin with 
\begin{lemma} \label{incl-1}
We have $\h_{nk} \subseteq \hs$.
\end{lemma}
\begin{proof}
Since the connection $\bnabla$ preserves both $g$ and $J$ we have that 
$\h_{nk}$ is a subalgebra of $\mathfrak{u}(\mathbb{V})$. Differentiating 
$\bnabla R=0$ (see (iii) in Lemma \ref{curv2}) shows, after using the Ricci identity that $R$ is invariant under 
the Lie algebra action of $\h_{nk}$. Explicitly 
$$\begin{array}{ll}
&R(FU_1,U_2,U_3,U_4)+R(U_1,FU_2,U_3,U_4)+R(U_1,U_2,FU_3,U_4)+\\
&R(U_1,U_2,U_3,FU_4)=0
\end{array}$$
whenever $F \in \h_{nk}$. Because $F$ is skew-symmetric with respect to $g$ this is equivalent to 
$$[F,R(U_1,U_2)]=R(FU_1,U_2)+R(U_1,FU_2).$$
In other words  
$$ [\h_{nk},\hs] \subseteq \hs 
$$
within $\mathfrak{u}(\mathbb{V})$. The claim follows now by a standard 
argument for compact Lie algebras. Consider the splitting 
$\mathfrak{u}(\mathbb{V})=\hs\oplus \hs^{\perp}$ with respect to the Killing form of the former. Since this is $\hs$-invariant we have $\h_{nk} \subseteq 
\hs \oplus \{F \in \hs^{\perp} : [\hs,F]=0 \}$. Because $(\hs,\mathbb{V})$ is irreducible and $J \in \hs$ the last summand vanishes by Schur's lemma.
\end{proof}
It follows that $\g_{nk} \subseteq \gs $, however this is {\it{not}} a Lie algebra 
inclusion since the Lie brackets on $\g_{nk}$ and $\gs $ 
are different. This situation requires some further identifications, as follows. 
Whenever $x \in \H$ let $\overline{x} \in \mathfrak{u}(\mathbb{V})$ be given by 
$$\overline{x}(v)=A_xv, \ \overline{x}(y)=A_xy.
$$
Note that $\overline{x}$ is orthogonal to $\mathfrak{u}(\V) \oplus \mathfrak{u}(\H)$.
\begin{lemma} \label{incl-2}
We have a well defined and injective map $\H \to \hs$ given by $x \in \H \mapsto \overline{x}$. Moreover, the induced map $\H \to \hs \slash \h_{nk}$ is injective.
\end{lemma}
\begin{proof}
By Corollary \ref{curv-s} we have $R(v,y)=-\overline{A_yv}$.
The equality $$\H=span \{A_yv:y \in \H, v \in \V\},$$  which is guaranteed by Lemma \ref{gammairr}, ensures the whole of $span \{\overline{x}, x \in \H \}$ is contained in $\hs$. To check both injectivity claims consider $x \in\H$ such that 
$\overline{x} \in \h_{nk}$. Because $\h_{nk}$ is contained in $\mathfrak{u}(\V) \oplus \mathfrak{u}(\H)$, whilst $\overline{x} \in \V \otimes H$, it follows that $\overline{x}=0$. The fact that $x=0$ follows again by Lemma \ref{gammairr}, (i). \\
\end{proof}
Let $\sigma: \g_{nk} \to \gs$ be given by 
$$ \sigma_{\vert \h_{nk} \oplus \V}=1_{\h_{nk} \oplus \V}, \ \sigma(x)=x-\overline{x}, x \in \H.
$$
This map is well-defined by (i) in Lemma \ref{incl-2} and has the following key properties.
\begin{propn} \label{incl-3}
We have that $\sigma$ is an injective Lie algebra morphism. Moreover:
\begin{itemize}
\item[(i)] $\gs=\hs+\sigma(\g_{nk})$
\item[(ii)] $\sigma(\g_{nk}) \neq \gs$.
\end{itemize}
\end{propn}
\begin{proof}
We compute 
\begin{align*}
[\sigma(x), \sigma(y)]_{sym}&=[x-\overline{x},y-\overline{y}] \\
&=  [x,y]_{sym}-[\overline{x},y]-[x,\overline{y}]+[\overline{x}, \overline{y}]\\
&= R(x,y)-2A_xy+[\overline{x}, \overline{y}]\\ 
\end{align*}

by using the various definitions. Since $[\overline{x},\overline{y}]=[A_x,A_y]$, using the comparison formula \eqref{curv11} leads to 
$$ [\sigma(x), \sigma(y)]_{sym}=\overline{R}(x,y)-2A_xy=[x,y]_{nk}=\sigma([x,y]_{nk}).
$$
Moreover 
$$\begin{array}{ll}
&[\sigma(v), \sigma(x)]_{sym}=[v,x-\overline{x}]_{sym}=R(v,x)+A_xv=
-\overline{A_xv}+A_xv.
\end{array}
$$
At the same time $\overline{R}$ is symmetric in pairs and thus $\overline{R}(v,x)=0$.
Therefore $$\sigma ([v,x]_{nk})=\sigma(A_xv)=A_xv-\overline{A_xv}$$ showing that $\sigma$ is a Lie algebra morphism on pairs of the type $(v,x)$. 

Choose now $F \in \h_{nk}$. Since $A$ is $\bnabla$-parallel, it is 
$\h_{nk}$-invariant, i.e. $[F,A_x]v=A_{Fx}v$ for all $x \in \H, v \in \V$. Equivalently, $[F,\overline{x}]=\overline{Fx}$ which leads easily to 
$[\sigma(F),\sigma(x)]_{sym}=\sigma([F,x]_{nk})$. Since $\sigma$ is the identity on 
$\h_{nk} \oplus \V$ we have finished showing that $\sigma$ is a Lie algebra morphism. Clearly $\sigma$ is injective so there remains only to prove (i)and (ii).\\
For (i), we have $\g_{sym}=\hs \oplus \mathbb{V}$. Since 
$\sigma(\g_{nk})=\h_{nk} \oplus \sigma(\mathbb{V})$ and $\h_{nk} \subseteq \hs$ 
we have $\hs +\sigma(\g_{nk})=\hs+\sigma(\mathbb{V})$. From the construction of 
$\sigma$ we know that $\mathbb{V} \subseteq \hs +\sigma(\mathbb{V})$ since 
$x=\sigma(x)+\overline{x}$ for all $x \in \H$ and $\sigma_{\vert \V}=1_{\V}$. This proves the claim.\\
Finally for (ii), assume equality holds. Since $\h_{nk} \subseteq  \sigma(\g_{nk})$ and $\hs 
\subseteq \gs$ have the same codimension and $\h_{nk} \subseteq \hs$ it follows that $\h_{nk}=\hs$. Then the last part of Lemma \ref{incl-2} yields $\H=\{0\}$, a contradiction.
\end{proof}
Note that in the proposition above the isotropy sub-algebra of $\g_{nk}$ is recovered from 
\begin{equation} \label{recov-i}
\h_{nk}=\hs \cap \sigma(\g_{nk}).
\end{equation}
In order to progress towards the last step of the argument we make the following observations. Since $(M,g)$ is a Hermitian locally symmetric space with $g$ locally irreducible, the pair of Lie algebras $(\gs,\hs)$ is an irreducible Hermitian symmetric pair of compact type. By the general theory of such pairs, $\gs$ is a simple compact Lie algebra and so its Killing form is negative definite. 

Identifying $\g_{nk}$ with a subalgebra of the compact Lie algebra $\gs$ makes it straightforward to show:
\begin{lemma} \label{semi-easy}
$\g_{nk}$ is a semisimple Lie algebra of compact type.  
\end{lemma}
\begin{proof}
Since $\g_{nk}$ is a subalgebra of a (semi-)simple Lie algebra of compact 
type it is semisimple if and only if its 
center $\mathfrak{z}(\g_{nk})$ vanishes identically. Moreover if that happens $\g_{nk}$ has also compact type. Therefore, pick $F+U \in \mathfrak{z}(\g_{nk})$ where $F \in \h_{nk}$ and 
$U \in \mathfrak{h}_{nk}$. From $[F+U,\h_{nk}]=0$ we get $[F,\h_{nk}]=0$ and $\h_{nk}U=0$. Because $(\h_{nk},\V)$ and $(\h_{nk},\H)$ are irreducible 
it follows that $U=0$. The vanishing of $F$ follows now from $[F,\mathbb{V}]=0$, thus proving the claim.
\end{proof}
%\begin{rem} \label{not-needed}
%It is possible, even at this stage, to show directly that $\g_{nk}$ is a simple 
%Lie algebra. However this is not necessary since the classification will %automatically imply this fact.
%\end{rem}
We now observe that $\mathbb{V}$ is intrinsically determined within $\g_{nk}$ as the orthogonal complement of $\h_{nk}$, with respect to the Killing form of $\g_{nk}$. 
This is standard theory for symmetric spaces; however the fact does not always hold in the case of arbitray  infinitesimal models coming from connections with torsion.
Proposition \ref{uniqueV} below will be used to determine $\mathbb{V}$ when the pair $(\g_{nk},\h_{nk})$ is already given explicitly.
\begin{propn} \label{uniqueV}
The representation $(\h_{nk},\mathbb{V})$ coincides with the isotropy representation of $\h_{nk}$ in $\g_{nk}$.
\end{propn}
\begin{proof}
We must show that $\bbV$ is the orthogonal complement 
$\h_{nk}^{\perp}$ of $\h_{nk}$ in $\g_{nk}$, taken with respect to the Killing form 
$Q$ of $\g_{nk}$.
Let $F \in \h_{nk}$ and $U \in \mathbb{V}$. Recall that $\h_{nk}$ is a subalgebra 
in $\g_{nk}$ preserving the splitting $\mathbb{V}=\V \oplus \H$. Also recall that 
the torsion satisfies $\tau(\V,\V)=0, \tau(\V,\H)\subseteq \H, \tau(\H, \H)\subseteq \V$. Based on these facts it follows directly from the definition of the Lie bracket $[\cdot, \cdot]_{nk}$ that 
\begin{equation*}
\begin{split}
&(ad_U \circ ad_F) \h_{nk} \subseteq \mathbb{V}, \ (ad_U \circ ad_F)\V \subseteq \h_{nk} \oplus \H \\
& (ad_U \circ ad_F)\H \subseteq \h_{nk} \oplus \V \ \mbox{if} \ U \in \H \\
&(ad_U \circ ad_F)\H \subseteq \H \ \mbox{if} \ U \in \V.
\end{split}
\end{equation*}
These relations lead easily to $trace(ad_U \circ ad_F)=0$ when $U \in \H$ respectively to 
$$trace(ad_U \circ ad_F)=\sum \limits_{i}g(\tau(U,Fe_i),e_i)
$$
when $U \in \V$, where $\{e_i\}$ is some $g$-orthonormal basis in $\H$. But $g(\tau(U,Fe_i),e_i)=
g(A_{Fe_i}U,e_i)=g(U,A_{e_i}Fe_i)$ and the contraction $\sum \limits_{i}A_{e_i}Fe_i=0$ since the tensor \newline  $A:\H \times \H \to \V$ is $J$-anti-invariant, whereas $FJ=JF$.  Thus, we have showed that $Q(\h_{nk},\mathbb{V})=0$ so it follows $\mathbb{V} \subseteq \h_{nk}^{\perp}$. The converse inclusion is an entirely standard argument which we only outline for completeness. Pick an element $F+U \in \h_{nk}^{\perp}$ where $F \in \h_{nk}, U \in \mathbb{V}$.
Because $\h_{nk}$ and $\mathbb{V}$ are $Q$-orthogonal it follows that 
$Q(F,\h_{nk})=0$. Since $Q(F,\mathbb{V})=0$ (again since $F \in \h_{nk}$) it follows that $F$ is $Q$-orthogonal 
to $\g_{nk}$. As $\g_{nk}$ is semisimple by Lemma \ref{semi-easy}, $F$ must vanish and the claim follows.
\end{proof}
With these in hand we can prove the main result in this section. For the convenience of the reader we recall in Table \ref{table1} below the classification of irreducible Hermitian symmetric pairs $(\gs,\hs)$ of Lie algebras of compact type.
 \begin{table}[h!]
\centering
\begin{tabular}{||c | c | c||} 
 \hline
$ \gs$ & $\hs$ & Restrictions \\ 
 [0.5ex] 
 \hline\hline
 $\mathfrak{su}(n)$  & $s(\mathfrak{u}(p)\oplus \mathfrak{u}(q))$ & $n\geq 2, \  p+q = n $\\ 
 $\mathfrak{so}(n+2)$ & $\mathfrak{so}(n)\oplus \mathfrak{so}(2)$ & $n\geq 2$\\
 $\mathfrak{sp}(n)$ & $\mathfrak{u}(n)$  & $n \geq 3$  \\
%$SO_{2n}/\mathbb{Z}_2$ & $(SO_{2(n-1)}\times SO_2) / \mathbb{Z}_2 $ & $n\geq 4$ \\
$\mathfrak{so}(2n)$ & $\mathfrak{u}(n)$ & $n\geq 4$ \\
 $\mathfrak{e}_6$ & $ \mathfrak{so}(10) \oplus \mathfrak{so}(2)$ & -- \\  
 $\mathfrak{e}_7$ & $\mathfrak{e}_6 \oplus \mathfrak{so}(2) $ & -- \\  [1ex]
    \hline
\end{tabular}
\vspace{2mm}
\caption{Irreducible Hermitian symmetric pairs of compact type}
\label{table1}
\end{table}
\bigskip

We can now prove the following:
\begin{thm} \label{direct}
We have either 
$$ \gs=\mathfrak{su}(2n), \hs=s(\mathfrak{u}(1)\oplus \mathfrak{u}(2n-1)), \ \g_{nk}=\mathfrak{sp}(n), \ \h_{nk}=\mathfrak{sp}(n-1) \oplus \mathfrak{u}(1)
$$
or 
$$ \gs=\mathfrak{so}(2n+2), \hs=\mathfrak{u}(n+1), \ \g_{nk}=\mathfrak{so}(2n+1), \ \h_{nk}=\mathfrak{u}(n).
$$
Moreover both cases do occur.
\end{thm}
\begin{proof}
In \cite{on} Onishchik gives the full classification for triples $(\g, \g^{\prime}, 
\g^{\prime \prime})$ where $\g$ is a simple compact Lie algebra and $\g^{\prime}, 
\g^{\prime \prime}$ are {\it{strict}} subalgebras of $\g$ such that $\g=\g^{\prime}+\g^{\prime \prime}$. By Proposition \ref{incl-3} the triple $(\gs, \hs, \sigma(\g_{nk}))$ satisfies exactly these conditions. We now take from the list in \cite{on} only the 
occurences when $(\g, \g^{\prime})$ is a Hermitian symmetric pair as listed above. Keeping in mind that $\h_{nk}$ is determined from \eqref{recov-i} we have the following possibilities.
 \begin{table}[h!]
\centering
\begin{tabular}{||c | c| c | c||} 
 \hline
 \hline
 $\gs$ & $\hs$ & $\g_{nk}$ & $\h_{nk}$ \\
 \hline
$\mathfrak{su}(2n)$ & $ s(\mathfrak{u}(1)\oplus \mathfrak{u}(2n-1))$ & $\mathfrak{sp}(n)$& $\mathfrak{sp}(n-1) \oplus \mathfrak{u}(1)$\\
$\mathfrak{so}(2n+2)$& $\mathfrak{u}(n+1)$ & $\mathfrak{so}(2n+1)$ & $\mathfrak{u}(n)$\\
$\mathfrak{so}(7)$& $\mathfrak{so}(2)\oplus \mathfrak{so}(5)$& $\mathfrak{g}_2$ & $\mathfrak{u}(2)$ \\
\hline
\end{tabular}
\caption{Onishchik's List}
\label{table3}
\end{table}

See also \cite{kerr}, page 155, for a nice account on how Onishchik's list compares 
to the list of symmetric spaces of compact type. For the last  entry in the table 
we split   $\mathfrak{g}_2=\mathfrak{u}(2) \oplus \mathfrak{m}_2$ with respect to the Killing forms of  $\mathfrak{g}_2$. By Proposition \ref{uniqueV} the representation  $(\mathfrak{u}(2),\mathfrak{m}_2)$  corresponds to the representation $(\h_{nk},\mathbb{V})$. However M.Kerr shows in \cite{kerr} that this  representation has $3$ inequivalent irreducible components which contradicts 
that $(\h_{nk},\mathbb{V})$ has only $2$ irreducible components, by Proposition \ref{irrep-fn} and Proposition \ref{irrep-V}. 

For the two remaining instances, we only need argue that the foliation is uniquely 
determined. This is because in both cases we already know that foliations of the type we require exist. Indeed, for $(\g_{nk},\h_{nk})=(\mathfrak{sp}(n), \mathfrak{sp}(n-1) \oplus \mathfrak{u}(1))$ or $(\g_{nk},\h_{nk})=(\mathfrak{so}(2n+1), \mathfrak{u}(n))$, consider the isotropy decomposition $\g_{nk}=\h_{nk} \oplus \mathfrak{m}$ with respect to the Killing form of $\g_{nk}$. Then (see e.g.\cite{kerr}) we have $\mathfrak{m}=\mathfrak{m}_1 \oplus 
\mathfrak{m_2}$ with $\mathfrak{m}_i, i=1,2$ irreducible, inequivalent 
representations of $\h_{nk}$. The space $\mathfrak{m}_1$ satisfies $[\mathfrak{m}_1, \mathfrak{m}_1]
_{\mathfrak{m}}=0$, while $[\mathfrak{m}_2, \mathfrak{m}_2]
_{\mathfrak{m}}=\mathfrak{m}_1$. \\ 
By Proposition \ref{uniqueV} we have $\bbV=\mathfrak{m}$. Since 
both representations $(\h_{nk},\V)$ and $(\h_{nk},\H)$ are irreducible we must 
have either $\V=\mathfrak{m}_1, \H=\mathfrak{m}_2$ or $\V=\mathfrak{m}_2, \H=\mathfrak{m}_1$. As the component on $\bbV$ of $[\V,\V]_{nk}$ vanishes(the torsion tensor vanishes on $\V$) only the first instance may occur. It corresponds to the Hopf fibration, respectively the twistor fibration of $\mathbb{S}^{2n}$ as explained  in \cite{kerr}. 

This determines completely the ambient space $N$ as well as the structure of the 
foliations, thus the proof of part (i) in Theorem \ref{t2} is now complete. 
\end{proof}

\endproof

\end{document}